\documentclass[10pt,a4paper]{article}%
\usepackage{amstext,amsfonts,amsmath,amssymb,makeidx}
\usepackage{bm}
\usepackage{graphicx}%
\usepackage[latin2]{inputenc}

\DeclareMathOperator{\sinc}{sinc}

\providecommand{\U}[1]{\protect\rule{.1in}{.1in}}
\newtheorem{theo}{Statement}
\newtheorem{theorem}{Theorem}

\newtheorem{conjecture}{Conjecture}
\newtheorem{corollary}{Corollary}

\newtheorem{remark}[theorem]{Remark}

\newenvironment{proof}[1][Proof]{\noindent\textbf{#1} }{\ \rule{0.5em}{0.5em}}

\newcommand\blfootnote[1]{%
  \begingroup
  \renewcommand\thefootnote{}\footnote{#1}%
  \addtocounter{footnote}{-1}%
  \endgroup
}

\begin{document}

\begin{center}

\centerline {\bf \large SOME IMPROVEMENTS OF JORDAN-STE\v CKIN}

\smallskip

\centerline {\bf \large  AND BECKER-STARK INEQUALITIES}
\bigskip
Marija Nenezi\' c${}^{\,\mbox{\tiny 1)}}$,
Ling Zhu${}^{\mbox{\tiny 2)}}$

\blfootnote{$\!\!\!\!\!\!\!\!\!\!\!{}^{*}\, $Corresponding author
}\blfootnote{$\!\!\!\!\!\!\!\!\!\!\!$E-mails: Marija~Nenezi\' c$\,<${\sl maria.nenezic@gmail.com}$>$,
Ling Zhu$\,<${\sl zhuling0571@163.com}$>$
}
\begin{center}
{\footnotesize \it
${}^{1)}$School of Electrical Engineering, University of Belgrade,                           \\[+0.20 ex]
Bulevar Kralja Aleksandra 73, 11000 Belgrade, Serbia                                          \\[+1.00 ex]
${}^{2)}$Department of Mathematics, Zhejiang Gongshang University,                            \\[+0.20 ex]
Hangzhou 310018, China                                                                        \\[+1.00 ex]
}
\end{center}
\[
\]

\end{center}

\bigskip
\noindent
{\small {\bf Abstract.}
The aim of this article is to give some improvements of Jordan-Ste\v ckin and Becker-Stark inequalities discussed in [1].
}

\bigskip
\noindent
{\footnotesize {\bf MSC:}  26D15, 41A10, 42A16 }

\bigskip
\noindent
{\footnotesize {\bf Keyword.}  Jordan-Ste\v ckin inequalities; Becker-Stark inequalities }

\section{Introduction}

\quad  L. \textsc{Debnath}, C. \textsc{Mortici} and L. \textsc{Zhu} in \cite{refinements} discussed about \textsc{Jordan}'s inequality:%
\begin{equation}
\frac{\sin x}{x} \geq \frac{2}{\pi},  \quad  x \in ( 0, \pi/2]
\end{equation}
and it's improvements

\begin{equation}
\frac{2}{\pi} + \frac{1}{\pi^{3}}\left(\pi ^ {2} - 4 x ^ 2 \right) \leq \frac{\sin x}{x} \leq \frac{2}{\pi} + \frac{\pi - 2}{\pi^{3}}\left(\pi ^ {2} - 4 x ^ 2 \right) ,  \quad  x \in ( 0, \pi/2],
\end{equation}
and

\begin{equation}
\frac{2}{\pi} + \frac{1}{2 \pi^{5}}\left(\pi ^ {4} - 16 x ^ 4 \right) \leq \frac{\sin x}{x} \leq \frac{2}{\pi} + \frac{\pi - 2}{\pi^{5}}\left(\pi ^ {4} - 16 x ^ 4 \right) ,  \quad  x \in ( 0, \pi/2].
\end{equation}

They concluded that equalities in (2) and (3) hold iff $x = \pi/2$. In the case when $x \to 0_+$, we have the equalities in the right - hand side of (2) and (3), and strict inequalities on the left - hand side of (2) and (3).

In \cite{refinements} (\textit{Theorem 1}, \textit{Theorem 2}) the left - hand side of (2) and (3) near zero were improved.

The following inequality:

\begin{equation}
\label{tan4}
\tan x \geq \frac{4}{\pi} \cdot \frac{x}{\pi - 2x},  \quad  x \in [0, \pi/2).
\end{equation}
well known as \textsc{Ste\v ckin}'s inequality was also analysed in \cite{refinements}.

As noted in \cite{refinements} this inequality becomes equality for $x = 0$, and
\begin{align*}
\lim_{x \to (\pi/2)_-} \left( \tan x - \frac{4}{\pi} \cdot \frac{x}{\pi-2x} \right) = \frac{2}{\pi}.
\end{align*}

Some improvements of (\ref{tan4}), in the left neighbourhood of $\pi/2$ were presented in \cite{refinements} (\textit{Theorem 3}, \textit{Theorem 4}).

\textsc{M. Becker} and \textsc{L. E. Stark} in \cite{becker} presented the inequality

\begin{equation}
\frac{8}{\pi ^ {2} - 4 x ^ {2}} < \frac{\tan x}{x} < \frac{\pi ^ {2}}{\pi ^ {2} - 4 x ^ {2}}, \quad 0 < x < \frac{\pi}{2}.
\end{equation}

Some double inequalities of the \textsc{Becker-Stark} type, were proposed in \cite{refinements} (\textit{Theorem 5}, \textit{Theorem 6}).

\smallskip
In this paper, we give generalizations and improvements of the inequalities stated in \textit{Theorem 1}, \textit{Theorem 2},
\textit{Theorem 3}, \textit{Theorem 4}, \textit{Theorem 5} and \textit{Theorem 6} from \cite{refinements}. They are cited below
for readers convenience.

\begin{theo}
$\left(\emph{\cite{refinements}}, \textit{Theorem 1} \right)$
For every $x \in (0,\pi/2)$, it holds
\begin{align}
\label{1}
& \frac{2}{\pi} + \frac{1}{\pi^{3}}\left(\pi ^ {2} - 4 x ^ 2 \right) + \left(1-\frac{3}{\pi} \right) - \left(\frac{1}{6}-\frac{4}{\pi ^ 3} \right) x ^ 2 < \notag \\
& < \frac{\sin x}{x} < \\
& < \frac{2}{\pi} + \frac{1}{\pi^{3}}\left(\pi ^ {2} - 4 x ^ 2 \right) + \left(1-\frac{3}{\pi} \right) - \left(\frac{1}{6}-\frac{4}{\pi ^ 3} \right) x ^ 2 + \frac{1}{120} x ^ 4. \notag
\end{align}
\end{theo}

\begin{theo}

$\left(\emph{\cite{refinements}}, \textit{Theorem 2} \right)$
For every $x \in (0,\pi/2)$, it holds
\begin{align}
\label{2}
& \frac{2}{\pi} + \frac{1}{2\pi^{5}}\left(\pi ^ {4} -  16 x ^ 4 \right) + \left(1-\frac{5}{2\pi} \right) - \frac{1}{6} x ^ 2 < \notag \\
& < \frac{\sin x}{x} < \\
& < \frac{2}{\pi} + \frac{\pi - 2}{\pi^{5}}\left(\pi ^ {4} -  16 x ^ 4 \right) + \left(1-\frac{5}{2\pi} \right) - \frac{1}{6} x ^ 2 + \left( \frac{8}{\pi ^ {5}} + \frac{1}{120}\right) x ^ {4}. \notag
\end{align}
\end{theo}

\begin{theo}
$\left(\emph{\cite{refinements}}, \textit{Theorem 3} \right)$
For every $x \in (0,\pi/2)$, it holds
\begin{align}
\label{3}
& \frac{2}{\pi} - \frac{1}{2} \left( \frac{\pi}{2} - x \right) < \tan x - \frac{4}{\pi} \cdot \frac {x}{\pi - 2x} < \frac{2}{\pi} - \frac{1}{3} \left( \frac{\pi}{2} - x \right).
\end{align}
\end{theo}

\begin{theo}
$\left(\emph{\cite{refinements}}, \textit{Theorem 4} \right)$
For every $x \in (0,1)$, it holds
\begin{align}
\label{4}
\left(1-\frac{4}{\pi ^ 2}\right) x - \frac{8}{\pi ^ 3} x ^ 2 < \tan x - \frac{4}{\pi} \cdot \frac{x}{\pi - 2x} < \left(1 - \frac{4}{\pi ^ 2} x\right).
\end{align}
\end{theo}

\begin{theo}
$\left(\emph{\cite{refinements}}, \textit{Theorem 5} \right)$
For every $x \in \left(0.373, \pi/2\right)$ in the left-hand side and for every $x \in (0.301, \pi/2)$ in the right-hand side, the following inequalities hold true:
\begin{align}
\label{5}
\frac{8 + a(x)}{\pi ^ {2} - 4 x ^ {2}} < \frac{\tan x}{x} < \frac{8 + b(x)}{\pi ^ {2} - 4 x ^ {2}},
\end{align}
where
\begin{align*}
a(x) = \frac{8}{\pi} \left(\frac{\pi}{2} - x\right) + \left(\frac{16}{\pi ^ {2}} - \frac{8}{3}\right)\left(\frac{\pi}{2} - x\right)^{2}
\end{align*}
and
\begin{align*}
b(x) = a(x) +  \left(\frac{32}{\pi ^ {3}} - \frac{8}{3\pi}\right)\left(\frac{\pi}{2} - x\right)^{3}.
\end{align*}

\end{theo}

\begin{theo}
$\left(\emph{\cite{refinements}}, \textit{Theorem 6} \right)$
For every real number $x \in \left(0, 1.371\right)$, the following inequality holds true:
\begin{align}
\label{6}
\frac{\tan x}{x} < \frac{\pi ^ {2} - \left(4 - \frac{1}{3}\pi ^ {2}\right)x ^ {2} - (\frac{4}{3} - \frac{2}{15}\pi ^ {2})x ^ {4}}{\pi ^ {2} - 4 x ^ {2}}\,.
\end{align}
\end{theo}

\section{Preliminaries}

\quad Let $T_n^{\varphi, a}(x)$ denotes Taylor polynomial of order $n \in N$, associated with the function $\varphi(x)$ at the point $x = a$. $\overline{T}_n^{\varphi, a}(x)$ and $\underline{T}_n^{\varphi, a}(x)$ represent the Taylor polynomial of order $n \in N$, associated with the function $\varphi(x)$ at the point $x = a$, in the case $T_n^{\varphi, a}(x)\geq \varphi (x)$, respectively $T_n^{\varphi, a}(x)\leq \varphi (x)$, for every $x\in(a,b)$. We call $\overline{T}_n^{\varphi, a}(x)$ and $\underline{T}_n^{\varphi, a}(x)$ an upward and a downward approximation of $\varphi$ on $(a,b)$, respectively.

\smallskip
As discused in paper \cite{maria} for the sine function the following inqualities hold:
\begin{equation}
\label{taylorsin}
\begin{array}{c}
\underline{T}_3^{\sin, 0}(x)
<
\underline{T}_7^{\sin, 0}(x)
<
\underline{T}_{11}^{\sin, 0}(x)
<
\underline{T}_{15}^{\sin, 0}(x)
<
\ldots
<
\sin x
<
\ldots \\[1.0 ex]
\;\;\;
<
\overline{T}_{13}^{\sin, 0}(x)
<
\overline{T}_9^{\sin, 0}(x)
<
\overline{T}_5^{\sin, 0}(x)
<
\overline{T}_1^{\sin, 0}(x),
\end{array}
\end{equation}
for $x \in (0, \sqrt{12}) = (0, 3.464...)$.

\smallskip
We have the following \textsc{Taylor} series of $\sinc x $:
\begin{equation}
\label{taylorsinc}
\sinc x = \frac{\sin x}{x} = \displaystyle\sum_{k = 0}^{\infty}(-1) ^ {k}\frac{x^{2k}}{(2k + 1)!}
\end{equation}
for $x \neq 0$.

\smallskip
According to \cite{ryzhik} for $x \in (0, \pi/2)$ we have the following series representations:
\begin{equation}
\label{seriestan}
\tan x = \sum_{k = 1}^{\infty}\frac{2^{2k}(2^{2k}-1)}{(2k)!}\mid \!\!B_{2k}\!\!\mid x^{2k-1}
\end{equation}
and

\begin{equation}
\label{seriescotan}
\cot x = \frac{1}{x} -\sum_{k = 1}^{\infty}\frac{2^{2k}\mid \!\!B_{2k}\!\!\mid}{(2k)!} x ^ {2k-1}
\end{equation}
where $B_i$ $(i \in N)$ are \textsc{Bernoulli}'s numbers.

\smallskip
Suppose that $f(x)$ is a real function on $(a,b)$, and that $n$ is a positive integer such that $f ^ {(k)} (a+)$, $f ^ {(k)} (b-)$,
$(k \in {0, 1, 2, \ldots , n-1})$ exist. Let us denote by:
\begin{align*}
\boldmath{\boldmath{\underline{\cal T}}}_n^{f;b,a}(x) & = \sum_{k = 0}^{n-1}\frac{f^{(k)}(b-)}{k!}(x - b) ^ k + \notag\\
& + \frac{1}{(a - b) ^ n}\left(f(a+)-\sum_{k = 0}^{n-1}\frac{(a-b) ^ k f ^ {(k)}(b-)}{k!}\right)(x - b) ^ n
\end{align*}
and
\begin{align*}
\boldmath{\boldmath{\overline{\cal T}}}_n^{f;a,b}(x) & = \sum_{k = 0}^{n-1}\frac{f ^ {(k)}(a+)}{k!}(x - a) ^ k  + \notag\\
& + \frac{1}{(b - a) ^ n} \left(f(b-)-\sum_{k = 0}^{n-1}\frac{(b - a) ^ k f ^ {(k)} (a+)}{k!}\right)(x - a) ^ n.
\end{align*}

\textsc{S. Wu} and \textsc{L. Debnath} in \cite{wu} proved the following theorem:

\bigskip

\textbf{Theorem WD}
\textit{
Suppose that $f(x)$ is a real function on $(a,b)$, and that $n$ is a positive integer such that $f ^ {(k)} (a+)$, $f ^ {(k)} (b-)$,
$(k \in {0, 1, 2, \ldots , n})$ exist.
\renewcommand{\labelenumi}{\textit{(\roman{enumi})}}
\begin{enumerate}
\item
Supposing that $(-1) ^ {(n)} f ^ {(n)} (x)$ is increasing on $(a,b)$, then for all $x \in (a,b)$ the following inequality holds :
\begin{align}
\boldmath{\boldmath{\underline{\cal T}}}_n^{\phi;b,a}(x) < f(x) < \overline{T}_n^{f,b}(x)
\end{align}
Furthermore, if $(-1) ^ n f ^ {(n)}(x)$ is decreasing on $(a,b)$, then the reversed inequality of holds.
\item
Supposing that $f ^ {(n)}(x)$ is increasing on $(a,b)$, then for all $x \in (a,b)$ the following inequality holds:
\begin{align}
\boldmath{\boldmath{\overline{\cal T}}}_n^{\phi;a,b}(x) > f(x) > \underline{T}_n^{f,a}(x).
\end{align}
Furthermore, if $f ^ {(n)}(x)$ is decreasing on $(a,b)$, then the reversed inequality of holds.
\end{enumerate}
}
Some interesting applications of the previous theorem can be found in \cite{tanja}, \cite{Malesevic_Rasajski_Lutovac_2017c}, \cite{Milica_Makragic_2017} and \cite{22}.

\section{Main Results}

\subsection{Improvements of inequalities in Statement 1}

According to (\ref{taylorsin}), we can approximate $\sinc x$ function as it follows :
\begin{equation}
\label{tsinc}
\begin{array}{c}
\underline{T}_2^{sinc, 0}(x)
<
\underline{T}_6^{sinc, 0}(x)
<
\underline{T}_{10}^{sinc, 0}(x)
<
\underline{T}_{14}^{sinc, 0}(x)
<
\ldots
<
\sinc x
<
\ldots \\[1.0 ex]
\;\;\;\;\;
<
\overline{T}_{12}^{sinc, 0}(x)
<
\overline{T}_8^{sinc, 0}(x)
<
\overline{T}_4^{sinc, 0}(x)
<
\overline{T}_0^{sinc, 0}(x),
\end{array}
\end{equation}
for $x \in (0, \pi/2)\subset (0, \sqrt{12})$.

\smallskip
Based on approximation (\ref{tsinc}) we have the following theorem:
\begin{theorem}
For every $x \in (0,\pi/2)$ we have:
\begin{align}
\label{1}
& \underline{T}_2^{sinc, 0}(x) = \frac{2}{\pi} + \frac{1}{\pi^{3}}\left(\pi ^ {2} - 4 x ^ 2 \right) + \left(1-\frac{3}{\pi} \right) - \left(\frac{1}{6}-\frac{4}{\pi ^ 3} \right) x ^ 2 \leq  \notag\\
& \leq \underline{T}_{4k_1-2}^{sinc, 0}(x) < \sinc x < \overline{T}_{4k_2}^{sinc, 0}(x) \leq \frac{2}{\pi} + \frac{1}{\pi^{3}}\left(\pi ^ {2} - 4 x ^ 2 \right) +\\
& + \left(1-\frac{3}{\pi} \right) - \left(\frac{1}{6}-\frac{4}{\pi ^ 3} \right) x ^ 2 + \frac{1}{120} x ^ 4 = \overline{T}_4^{sinc, 0}(x) < \overline{T}_0^{sinc, 0}(x), \notag
\end{align}
for $k_1, k_2 \in N$.
\end{theorem}
\setcounter{theorem}{0}
\begin{remark}
It is obvious that Statement 1 is the special case of Theorem 1.
\end{remark}

\subsection{Improvements of inequalities in Statement 2}

Consider the following polynomials from inequality (\ref{2}) from \textit{Statement 2}:
\begin{align*}
Q_4(x) = \frac{2}{\pi} + \frac{1}{2\pi^{5}}\left(\pi ^ {4} - 16 x ^ 4 \right) + \left(1-\frac{5}{2\pi} \right) - \frac{1}{6} x ^ 2 = -\frac{8 x ^ 4}{\pi ^ 5} - \frac{x ^ 2}{6} + 1
\end{align*}
and
\begin{align*}
R_4(x) & = \frac{2}{\pi} + \frac{\pi - 2}{\pi^{5}}\left(\pi ^ {4} - 16 x ^ 4 \right) + \left(1-\frac{5}{2\pi} \right) - \frac{1}{6} x ^ 2 + \left(\frac{8}{\pi ^ 5} + \frac{1}{120}\right)x ^ 4 \\
& = \left(-\frac{16}{\pi ^ 4} + \frac{40}{\pi ^ 5} + \frac{1}{120}\right)x ^ 4 - \frac{x ^ 2}{6} - \frac{5}{2\pi} + 2.
\end{align*}
We have the following theorem:
\begin{theorem}
For every $x \in (0,\pi/2)$ we have:
\begin{align}
\label{theorem2}
Q_4(x)\!<\!\underline{T}_6^{sinc, 0}(x)\! \leq \!\underline{T}_{4k_1-2}^{sinc, 0}(x)\!<\!\sinc x \!<\!\overline{T}_{4k_2}^{sinc, 0}(x)\! \leq \!\overline{T}_4^{sinc, 0}(x)\!<\!R_4(x),
\end{align}
for $k_1, k_2 \in N$.
\end{theorem}

\bigskip

\begin{proof}
In order to prove (\ref{theorem2}) it is sufficient to prove for every $x \in (0,\pi/2)$ that inequalities
$Q_4(x) < \underline{T}_6^{sinc, 0}(x)$ and $\overline{T}_4^{sinc, 0}(x) < R_4(x)$ are true.

\medskip
According to (\ref{taylorsinc}) we have:
\begin{align*}
& \overline{T}_4^{sinc, 0}(x) = 1 - \frac{x ^ 2}{6} + \frac{x ^ 4}{120} , \\
& \underline{T}_6^{sinc, 0}(x) = 1 - \frac{x ^ 2}{6} + \frac{x ^ 4}{120} -  \frac{x ^ 6}{5040}.
\end{align*}
It is obvious that
\begin{align*}
\underline{T}_6^{sinc, 0}(x) - Q_4(x) & >  \left(1 - \frac{x ^ 2}{6} + \frac{x ^ 4}{120} -  \frac{x ^ 6}{5040}\right) - \left(-\frac{8 x ^ 4}{\pi ^ 5} - \frac{x ^ 2}{6} + 1\right) = \\
& = \left(\frac{1}{120} + \frac{8}{\pi ^ 5}\right) x ^ 4 - \frac{x ^ 6}{5040} > 0
\end{align*}
and
\begin{align*}
R_4(x) - \overline{T}_4^{sinc, 0}(x) & > \left(-\frac{16}{\pi ^ 4} + \frac{40}{\pi ^ 5} + \frac{1}{120}\right)x ^ 4 - \frac{x ^ 2}{6} - \frac{5}{2\pi} + 2 \\
& - \left(1 - \frac{x ^ 2}{6} + \frac{x ^ 4}{120}\right) = \left(-\frac{16}{\pi^4} + \frac{40}{\pi^5}\right) x ^ 4 - \frac{5}{2\pi} + 1 > 0
\end{align*}
hold for $x\in(0,\pi/2)$.
\end{proof}
\setcounter{theorem}{1}
\begin{remark}
Statement 2 is the special case of Theorem 2.
\end{remark}

\subsection{Improvements of inequalities in Statement 3}

In the monography \cite{mitrinovic}, \textsc{D.S.Mitrinovi\'c} discused about \textsc{Ste\v ckin's} inequality:
\begin{align*}
\tan x >  \frac{4}{\pi} \cdot \frac{x}{\pi - 2x},
\end{align*}
for $x \in (0,\pi/2)$.
Let us denote:
\begin{align}
\label{f3}
f(x) = \tan x - \frac{4 x}{\pi\left(\pi - 2 x \right)},
\end{align}
for $x \in (0,\pi/2)$ and let us notice:
\begin{align*}
\lim_{x \to \pi/2-} f(x) = \frac{\pi}{2}.
\end{align*}

In \cite{refinements} inequalities (\ref{3}) were proposed as adequate approximations of function $f(x)$ in left neighbourhood of the point $x = \pi/2$.

\smallskip
By replacing $x$ with $\pi/2 - t$ in the function $f(x)$, we obtain the following:
\begin{align*}
g(t) = f\left(\frac{\pi}{2} - t\right) = \cot t - \frac{1}{t} + \frac{2}{\pi},
\end{align*}
for $t \in (0, \pi/2)$. According to (\ref{seriescotan}) we have that
\begin{align*}
\cot t < \overline{T}_n ^ {\cot,0} (t) = \frac{1}{t} - \sum_{k = 1}^{n}\frac{2^{2k} \mid \!B_{2k}\!\mid}{(2k)!} t ^ {2k-1}
\end{align*}
for $t \in \left(0,\pi/2\right]$ and $n \in N$. Further, we have the following:
\begin{align}
\label{gt}
g(t) < \overline{T}^{\cot,0}_{n}(t) - \frac{1}{t} + \frac{2}{\pi}
\end{align}
and according to \textit{Theorem WD}
\begin{align}
\label{ct}
\mbox{\rm cot}\,t > \mbox{$\underline{\boldmath{\boldmath{\cal T}}}$}^{\cot;0,\pi/2}_{\!n}(t)
=
T^{\cot,0}_{n-1}(t) + \left(\mbox{\small $\displaystyle\frac{2}{\pi}$}\right)^{\!n}
\!{\bigg (}g\left(\mbox{\small$\displaystyle\frac{\pi}{2}$}\right) - T^{\cot,0}_{n-1}
\left(\mbox{\small $\displaystyle\frac{\pi}{2}$}\right)\!{\bigg )} t^{n},
\end{align}
for $t \in \left(0,\pi/2\right]$ and $n \in N$.
According to (\ref{gt}) and (\ref{ct}) we have the following:
\begin{align}
\label{gt2}
g(t) > \mbox{$\underline{\boldmath{\boldmath{\cal T}}}$}^{\cot;0,\pi/2}_{\!n}(t) - \frac{1}{t} + \frac{2}{\pi}
\end{align}
for $t \in \left(0,\pi/2\right]$.
Let us denote by:
\begin{align*}
F_n^{g}(t) = \overline{T}^{\cot,0}_{n}(t) - \frac{1}{t} + \frac{2}{\pi}
\end{align*}
and
\begin{align*}
{\cal F}_n^{g}(t) = \mbox{$\underline{\boldmath{\boldmath{\cal T}}}$}^{\cot;0,\pi/2}_{\!n}(t) - \frac{1}{t} + \frac{2}{\pi}.
\end{align*}
Returning replacement $t = \pi/2-x$ in (\ref{gt}) and (\ref{gt2}), we have the following theorem:
\begin{theorem}
For $x \in (0,\pi/2)$ and $n \in N$ we have:
\begin{align}
{\cal F}^{g}_{\!n}\!\left(\frac{\pi}{2} - x \right) < f(x) < F_n^{g}\!\left(\frac{\pi}{2} - x \right)
\end{align}
\end{theorem}
\begin{corollary}
We have the following improvements for inequality (\ref{3}) given in \textit{Statement 3}.
\begin{enumerate}
\item For $n = 1$ and for $x \in (0, \pi/2)$ we have:
\begin{align*}
& Q_1(x) < {\cal F}^{g}_{\!1}\!\left(\frac{\pi}{2} - x \right) = \notag \\
= \frac{2}{\pi} - \frac{4}{\pi ^ 2} &\left(\frac{\pi}{2} - x\right) < f(x) < \frac{2}{\pi} - \frac{1}{3} \left(\frac{\pi}{2} - x\right) = \\
& = F_1^{g}\!\left(\frac{\pi}{2} - x \right) = R_1(x)\notag.
\end{align*}
\item For $n = 3$ and for $x \in (0, \pi/2)$ we have:
\begin{gather*}
\!\!\!\!\!\!\!\!\!\!\!\!Q_1(x)\!<\!{\cal F}^{g}_{\!1}\!\left(\frac{\pi}{2}\!-\!x\right)\!<\!{\cal F}^{g}_{\!3}\!\left(\frac{\pi}{2}\!-\!x\right)\!=\!\frac{2}{\pi}\!-\!\frac{1}{3}\!\left(\frac{\pi}{2}\!-\!x\right)\!-\!\left(\frac{2}{\pi}\right)^3\!\!\left(\frac{2}{\pi}\!-\!\frac{\pi}{6}\right)\!\!\left(\frac{\pi}{2}\!-\!x\right)^3 \notag \\
\!\!\!\!\!\!\!\!\!<\!f(x)\!<\!\frac{2}{\pi}\!-\!\frac{1}{3}\left(\frac{\pi}{2}\!-\!x\right)\!-\!\frac{\left(\frac{\pi}{2}\!-\!x\right)}{45}\!=\!F_3^{g}\!\left(\frac{\pi}{2}\!-\!x \right)\!<\!F_1^{g}\!\left(\frac{\pi}{2}\!-\!x\!\right)\!=\!R_1(x) \notag.
\end{gather*}

\end{enumerate}
\end{corollary}

\bigskip

\subsection{Improvements of inequalities in Statement 4}

For the function $f(x)$ defined in (\ref{f3}) and according to Taylor series of $\tan x$ function in (\ref{seriestan}) and the binomial expansion of $\frac{1}{1-\left(\frac{2}{\pi}x\right)}$ over interval $(0, \pi/2)$ we have:
\begin{equation}
\begin{array}{rl}
f(x) \!\!&\!\! = \tan x - \mbox{\small $\displaystyle\frac{4}{\pi}$} \cdot \mbox{\small $\displaystyle\frac{x}{\pi-2x}$}                                                                                                                                                         \\[1.5 ex]
     \!\!&\!\! = \displaystyle\sum\limits_{i=1}^{\infty} \displaystyle\frac{2^{2i}\left(2^{2i}-1\right)|B_{2i}|}{(2i)!} x ^ {2i-1} - \displaystyle\displaystyle\frac{4}{\pi^2} \cdot \displaystyle\frac{x}{1-\left(\mbox{\footnotesize $\displaystyle\frac{2}{\pi}$} x\right)}   \\[1.5 ex]
     \!\!&\!\! = \displaystyle\sum\limits_{i=1}^{\infty}{\displaystyle\frac{2^{2i}\left(2^{2i}-1\right)|B_{2i}|}{(2i)!} x ^{2i-1}} - \displaystyle\sum\limits_{j=1}^{\infty}\displaystyle\displaystyle\frac{2^{j+1}}{\pi^{j+1}} x^{j}                                            \\[1.5 ex]
     \!\!&\!\! = \displaystyle\sum\limits_{k=1}^{\infty}{\alpha_k x^{k}},
\end{array}
\end{equation}
where
\begin{align*}
\alpha_{k}  = \left\{\begin{array}{ccc}
-\mbox{\small $\displaystyle\frac{2^{k+1}}{\pi^{k+1}}$} &\!\!:\!\!& k \!=\! 2\ell                                              \\[2.0 ex]
\mbox{\small $\displaystyle\frac{2^{2k+1}\left(2^{2k+1}-1\right)|B_{k+1}|}{(k+1)!}$}
-
\mbox{\small $\displaystyle\frac{2^{k+1}}{\pi^{k+1}}$} &\!\!:\!\!& k \ = \! 2 \ell-1
\end{array}
\right.
\end{align*}
for $\ell \in N$. It easy to notice that:
\begin{align}
\alpha_{k} \!<\! 0 \;\;\mbox{for}\;\; k\! = \!2\ell; \;\;\mbox{other}\;\; \alpha_{k} \!>\! 0,
\end{align}
for $k \in N$. It is not hard to check that $\displaystyle \lim_{k \to \infty} \alpha_k = 0$ and $\displaystyle (\alpha_k)\!\!\searrow\,$.

Finally, based on (26) and (27) and based on \textsc{Leibnitz} theorem, we have the following theorem:

\begin{theorem}
For every $x \!\in\! \left(0,1\right)$ and $\ell \!\in\! N$ holds:
\begin{equation}
\underline{T}^{f,0}_{2\ell}\!\left(x\right) < f(x) < \overline{T}^{f,0}_{2\ell-1}\!\left(x\right).
\end{equation}
\end{theorem}
\setcounter{theorem}{2}
\begin{remark}
Inequality (28) for $\ell = 1$ presents inequality (\ref{4}) from \textit {Statement 4}.
\end{remark}
\setcounter{theorem}{0}

\bigskip

\subsection{Improvements of inequalities in Statement 5}

Consider the following function:
$$
\varphi(x) = \left(\pi^2-4x^2\right)\frac{\tan x}{x},
$$
for $x \in \left(0,\pi/2\right)$.

\smallskip
By replacing $x$ with $\pi/2 - t$ in the function $\varphi(x)$, we obtain the following:
$$
\psi(t) = \varphi\!\left(\mbox{\small $\displaystyle\frac{\pi}{2}$}-t\right) = \frac{8\,t\,(\pi-t)\,\mbox{\rm cot}\, t}{\pi-2\,t}
$$
for $t \!\in\! \left(0,\pi/2\right)$. Improvement or inequalities from (\ref{5}) are given with the following theorem:

\setcounter{theorem}{4}

\begin{theorem}
For every $x \!\in\! \left(0,\pi/2\right)$ holds$:$
$$
\begin{array}{l}
\underline{T}_{4}^{\,\psi, 0}\!\left(\mbox{\small $\displaystyle\frac{\pi}{2}$}\!-\!x\right)
=
8
+
\mbox{\small $\displaystyle\frac{8}{\pi}$}\!\left(\mbox{\small $\displaystyle\frac{\pi}{2}$}\!-\!x\right)
+
\left(
\mbox{\small $\displaystyle\frac{16}{\pi^{2}}$}
-
\mbox{\small $\displaystyle\frac{8}{3}$}
\right) \! \left(\mbox{\small $\displaystyle\frac{\pi}{2}$}\!-\!x\right)^{\!2}
+
\left(
\mbox{\small $\displaystyle\frac{32}{\pi^{3}}$}
-
\mbox{\small $\displaystyle\frac{8}{3\pi}$}
\right) \! \left(\mbox{\small $\displaystyle\frac{\pi}{2}$}\!-\!x\right)^{\!3}   \\[2.0 ex]
+
\left(
\mbox{\small $\displaystyle\frac{64}{\pi^{4}}$}
\!-\!
\mbox{\small $\displaystyle\frac{16}{3\pi^{2}}$}
\!-\!
\mbox{\small $\displaystyle\frac{8}{45}$}
\right) \! \left(\mbox{\small $\displaystyle\frac{\pi}{2}$}\!-\!x\right)^{\!4} < \\[3.0 ex]
<
\varphi(x)
<                                                                                                                               \\[1.5 ex]
<
\overline{T}_{5}^{\,\psi, 0}\!\left(\mbox{\small $\displaystyle\frac{\pi}{2}$}\!-\!x\right)
=
8
+
\mbox{\small $\displaystyle\frac{8}{\pi}$}\!\left(\mbox{\small $\displaystyle\frac{\pi}{2}$}\!-\!x\right)
+
\left(
\mbox{\small $\displaystyle\frac{16}{\pi^{2}}$}
\!-\!
\mbox{\small $\displaystyle\frac{8}{3}$}
\right) \! \left(\mbox{\small $\displaystyle\frac{\pi}{2}$}\!-\!x\right)^{\!2}
+
\left(
\mbox{\small $\displaystyle\frac{32}{\pi^{3}}$}
\!-\!
\mbox{\small $\displaystyle\frac{8}{3\pi}$}
\right) \! \left(\mbox{\small $\displaystyle\frac{\pi}{2}$}\!-\!x\right)^{\!3}   \\[2.0 ex]
+
\left(
\mbox{\small $\displaystyle\frac{64}{\pi^{4}}$}
\!-\!
\mbox{\small $\displaystyle\frac{16}{3\pi^{2}}$}
\!-\!
\mbox{\small $\displaystyle\frac{8}{45}$}
\right) \! \left(\mbox{\small $\displaystyle\frac{\pi}{2}$}\!-\!x\right)^{\!4}
+
\left(
\mbox{\small $\displaystyle\frac{128}{\pi^{5}}$}
\!-\!
\mbox{\small $\displaystyle\frac{32}{3\pi^{3}}$}
\!-\!
\mbox{\small $\displaystyle\frac{8}{45\pi}$}
\right) \! \left(\mbox{\small $\displaystyle\frac{\pi}{2}$}\!-\!x\right)^{\!5}.
\end{array}
$$
\end{theorem}
One proof of this statement is based on equivalent mixed trigonometric polynomial inequalities$:$
$$
f(x)
\!=\!
\left(\pi^2-4 x^2\right) \sin x
-
x \, T_{4}^{\,\psi, 0}\!\left(\mbox{\small $\displaystyle\frac{\pi}{2}$}\!-\!x\right) \cos x
> 0
$$
and
$$
g(x)
\!=\!
\left(\pi^2-4 x^2\right) \sin x
-
x \, T_{5}^{\,\psi, 0}\!\left(\mbox{\small $\displaystyle\frac{\pi}{2}$}\!-\!x\right) \cos x
< 0,
$$
for $x \!\in\! \left(0,\pi/2\right)$.
Papers \cite{Malesevic_Makragic_2016} and \cite{Lutovac_Malesevic_Mortici_2017} show that problem of proving mixed trigonometric polynomial inequalities is a deciable problem and these inequalities for the mixed trigonometric polynomial functions are followed by the algorithms from papers above. Some interesting applications of the algorithmic approach in proving mixed trigonometric inequalities can be found in \cite{21} and \cite{Malesevic_Lutovac_Banjac_2018}; see also \cite{Malesevic_Rasajski_Lutovac_2017a} and \cite{Malesevic_Rasajski_Lutovac_2017b}.
\setcounter{theorem}{3}
\begin{remark}
It is obvious that Statement 5 is consequence of Theorem 5.
\end{remark}
Further, let us observe array $(\alpha_k)_{k \in N}$ defined with:
$$
\alpha_1 \!=\! 1, \,
\alpha_{2j} \!=\! 0, \,
\alpha_{2j+1} \!=\! - \mbox{\footnotesize $\displaystyle\frac{2^{2j}|B_{2j}|}{(2j)!}$}
$$
for $j \!\in\! N$. Then based on \cite{ryzhik} we have the following series representations:
$$
\begin{array}{rcl}
\psi(t)
\!\!&\!\!=\!\!&\!\!
\mbox{\small $\displaystyle\frac{8}{\pi}$}\,t\,(\pi-t)
\,
\mbox{\small $\displaystyle\frac{1}{1-\left(\frac{2\,t}{\pi}\right)}$} \,\mbox{\rm cot}\, t                   \\[2.0 ex]
\!\!&\!\!=\!\!&\!\!
\mbox{\small $\displaystyle\frac{8}{\pi}$}\,t\,(\pi-t)
{\bigg (}\displaystyle\sum_{i=0}^{\infty}{\!\left(\mbox{\small $\displaystyle\frac{2\,t}{\pi}$}\right)^{\!i}}{\bigg )}\!
{\bigg (}\displaystyle\sum_{j=0}^{\infty}{\alpha_{2j+1}t^{2j-1}}{\bigg )}
\end{array}
$$
for $t \in \left(0,\pi/2\right)$.
Let $r_2(m)$ be the remainder after division of natural number $m$ with 2. We are posing the following conjecture:

\setcounter{theorem}{0}

\begin{conjecture}
\quad
\begin{enumerate}
\item For the function $\psi(t)$ on $t \!\in\! (0,\pi/2)$ the following equality holds$:$
\begin{equation}
\label{Conj_1}
\psi(t)
=
\displaystyle\sum_{m=0}^{\infty}{
{\Bigg (}
\frac{8\,\alpha_{m+1-r_2(m)}}{\pi^{r_2(m)}}
+\!
\displaystyle\sum_{i=1}^{[m/2]}{\frac{2^{2i+2+r_2(m)}\alpha_{m + 1 - 2i - r_2(m)}}{\pi^{2i+r_2(m)}}}\!{\Bigg )}t^{m}}.
\end{equation}

\item For the function $\psi(t)$ on $t \!\in\! (0,\pi/2)$ and $ \ell \in N$ the following inequalities are true:
\begin{equation}
\label{Conj_2}
\underline{T}_{2 \ell}^{\,\psi, 0}(t)
<
\psi(t)
<
\overline{T}_{2\ell + 1}^{\,\psi, 0}(t)
\quad
\left(\,t \!\in\! (0,\pi/2) \,\wedge\, \ell \!\in\! N\,\right).
\end{equation}
\end{enumerate}
\end{conjecture}

\bigskip

\subsection{Improvements of inequality in Statement 6}

Let us denote the following function:
\begin{align*}
f(x) = \left(\pi^2-4x^2\right)\frac{\tan x}{x}
\end{align*}
for $x \in (0, \pi/2)$.

\break

According to \cite{ryzhik} and (\ref{seriestan}) we have:
\begin{align}
f(x) = \sum_{k = 1}^{\infty}C_k x^{2k-2}
\end{align}
where
\begin{align}
C_k
=
\frac{\pi^2 \cdot 2^{2k}(2^{2k}-1)\mid\! B_{2k} \!\mid}{(2k)!}
-
\frac{4 \cdot 2^{2k-2}(2^{2k-2}-1)\mid\! B_{2k-2} \!\mid}{(2k-2)!},
\end{align}
and $x \in (0,c)$ and $0 < c < \pi/2$.
\setcounter{theorem}{5}
It is not hard to check $C_k < 0$ for $k \in N$.

\medskip
Finally, based on \textit{Theorem WD} we get the following theorem:
\begin{theorem}
For every $x\in(0,c)$, where $0 \!<\! c \!<\! \pi/2$, the following inequalities hold:
\begin{gather*}
\underline{{\cal T}}^{f;0,c}_{m_1}(x)
=\!\!\!
\sum_{k = 1}^{m_1-1}\!\!C_k x^{2k\!-\!2}
+
\left(\frac{1}{c}\right)^{\!\!2m_1\!-\!2}\!\!\left(f(c)
-
\left(\sum_{k = 1}^{m_1\!-\!1}C_k c^{2k\!-\!2}\right)\!\!\right)\!x^{2m_1\!-\!2} \\
< f(x) <\sum_{k = 1}^{m_2}C_k x^{2k-2} = \overline{\mbox{$T$}}^{f,0}_{m_2}\left(x\right),
\end{gather*}
for $m_1, m_2 \in N$.
\end{theorem}
\setcounter{theorem}{4}
\begin{remark}
It is obvious that Statement 6 is consequence of Theorem 6 for $m_2 = 3$.
\end{remark}
\qquad Approximations, discused in this paper, can have great significance for potential applications of analytic inequalities in engineering.
Some specific inequalities of the similar type are considered in \cite{Rahmatollahi_DeAbreu_2012}, \cite{Alirezaei_Mathar_2014} and \cite{Cloud_Drachman_Lebedev_2014}.

\bigskip

\textbf{Acknowledgements.}
We thank the anonymous reviewers for their careful reading of our manuscript and their many insightful comments and suggestions.
The second author was partially supported by the National Natural Science Foundation of China (no. 11471285 and no.  61772025).

%


\begin{thebibliography}{50}

\bibitem{refinements}
{\sc L. Debnath}, {\sc C. Mortici}, {\sc L. Zhu}: {\em Refinements of Jordan-Ste\v ckin and Becker-Stark inequalities},
Results Math. {\bf 67} (1-2), 207-215 (2015)

\bibitem{becker}
{\sc M. Becker}, {\sc E. L. Stark}: {\em On hierarchy of polynomial inequalities for tan(x)},
Univ. Beograd, Publ. Elektrotehn. Fak. Ser. Mat. Fiz. {\bf 602}-{\bf 633}, 133-138 (1978)

\bibitem{mitrinovic}
{\sc D.$\,$S. Mitrinovi\' c}: {\em Analytic inequalities}, Springer--Verlag (1970)

\bibitem{maria}
{\sc M. Nenezi\' c}, {\sc B. Male\v sevi\' c}, {\sc C. Mortici}: {\em New approximations of some expressions involving trigonometric functions}, Appl. Math. Comput. {\bf 283}, 299-315 (2016)

\bibitem{tanja}
{\sc B. Male\v sevi\' c}, {\sc T. Lutovac}, {\sc M. Ra\v sajski}, {\sc C. Mortici}:
{\em Extensions of the natural approach to refinements and generalizations of some trigonometric inequalities},  arXiv:1712.06792 (2017)

\bibitem{ryzhik}
{\sc I. Gradshteyn}, {\sc I. Ryzhik}: {\em Table of Integrals Series and Products}, 8-th Edition,
Academic Press (2015)

\bibitem{wu}
{\sc S. Wu}, {\sc L. Debnath}:
{\em A generalization of L'Hospital-type rules for monotonicity and its application}, Appl. Math. Lett. {\bf 22}, 284-290 (2009)


\bibitem{Zhu_2006}
{\sc L. Zhu}: {\em Sharpening of Jordan's inequalities and its applications}, Math. Inequal. Appl. {\bf 9} (1), 103-106 (2006)

\bibitem{Zhu_2006_II}
{\sc L. Zhu}: {\em Sharpening Jordan's inequality and Yang Le's inequality II}, Appl. Math. Lett. {\bf 19} (9), 990-994 (2006)


\bibitem{Mortici_2011}
{\sc C. Mortici}: {\em The natural approach of Wilker-Cusa-Huygens inequalities}, Math. Inequal. Appl. {\bf 14}:3, 535-541 (2011)

\bibitem{Rahmatollahi_DeAbreu_2012}
{\sc G. Rahmatollahi}, {\sc G.T.F. De Abreu}:
{\em Closed-Form Hop-Count Distributions in Random Networks with Arbitrary Routing}, IEEE Trans. Commun. {\bf 60}:2, 429-444 (2012)

\bibitem{Alirezaei_Mathar_2014}
{\sc G. Alirezaei}, {\sc R. Mathar}:
{\em Scrutinizing the average error probability for nakagami fading channels} in The IEEE International Symposium on Information Theory (ISIT'14), Honolulu, Hawai, USA, Jun. 2884-2888 (2014)

\bibitem{Cloud_Drachman_Lebedev_2014}
{\sc M.$\,$J. Cloud, B.$\,$C. Drachman, L.$\,$P. Lebedev}: {\em Inequalities With Applications to Engineering}, Springer (2014)

\bibitem{Nishizawa_2015}
{\sc Y. Nishizawa}: {\em Sharpening of Jordan's type and Shafer-Fink's type inequalities with exponenti alapproximations}, Appl. Math. Comput. {\bf 269}, 146-154 (2015)

\bibitem{Malesevic_Makragic_2016}
{\sc B. Male\v sevi\' c}, {\sc M. Makragi\' c}: {\em A Method for Proving Some Inequalities on Mixed Trigonometric Polynomial Functions}, J. Math. Inequal. {\bf 10}:3, 849-876 (2016)


\bibitem{Lutovac_Malesevic_Mortici_2017}
{\sc T. Lutovac}, {\sc B. Male\v sevi\' c}, {\sc C. Mortici}:
{\em The natural algorithmic approach of mixed trigonometric-polynomial problems},
J. Inequal. Appl. {\bf 2017}:116, 1-16 (2017)

\bibitem{Malesevic_Rasajski_Lutovac_2017a}
{\sc B. Male\v sevi\' c}, {\sc M. Ra\v sajski}, {\sc T. Lutovac}:
{\em Refinements and generalizations of some inequalities of Shafer-Fink's type for the inverse sine function},
J. Inequal. Appl. {\bf 2017}:275, 1-9 (2017)

\bibitem{Malesevic_Rasajski_Lutovac_2017b}
{\sc B. Male\v sevi\' c}, {\sc M. Ra\v sajski}, {\sc T. Lutovac}:
{\em Refined estimates and generalizations of inequalities related to the arctangent function and Shafer's inequality},
arXiv:1711.03786 (2017)

\bibitem{Malesevic_Rasajski_Lutovac_2017c}
{\sc B. Male\v sevi\' c}, {\sc M. Ra\v sajski}, {\sc T. Lutovac}:
{\em A new approach to the sharpening and generalizations of Shafer-Fink and Wilker type inequalities},
arXiv:1712.03772 (2017)

\bibitem{Milica_Makragic_2017}
{\sc M. Makragi\' c}: {\em A method for proving some inequalities on mixed hyper\-bolic-trigonometric polynomial functions},
J. Math. Inequal. {\bf 11}:3, 817-829 (2017)


\bibitem{Malesevic_Lutovac_Banjac_2018}
{\sc B. Male\v sevi\' c}, {\sc T. Lutovac}, {\sc B. Banjac}:
{\em A proof of an open problem of Yusuke Nishizawa for a power-exponential function},
Accepted in J. Math. Inequal. (2018)


\bibitem{3}
{\sc L. Debnath}, {\sc C.J. Zhao}: {\em New strengthened Jordan's inequality and its applications},
 Appl. Math. Lett. {\bf 16}, 557-560 (2003)



\bibitem{7}
{\sc C. Mortici}: {\em A subtly analysis of Wilker inequality},
 Appl. Math. Comput. {\bf 231}, 516-520 (2014)

\bibitem{8}
{\sc J. Pe\v cari\' c}, {\sc A.U. Rehman}: {\em Cauchy means introduced by an inequality of Levin and Steckin},
 East J. Approx. {\bf 15}, 515-524 (2009)



\bibitem{12}
{\sc Zh.-J. Sun}, {\sc L. Zhu}: {\em  Simple proofs of the Cusa-Huygens-type and Becker-Starktype inequalities},
 J. Math. Inequal. {\bf 7}(4), 563-567 (2013)


\bibitem{14}
{\sc S. Wu}, {\sc L. Debnath}: {\em A new generalized and sharp version of Jordan's inequality and its applications to the improvement of the Yang Le inequality},
 Appl. Math. Lett. {\bf 19}(12), 1378-1384 (2006)

\bibitem{15}
{\sc S. Wu}, {\sc L. Debnath}: {\em A new generalized and sharp version of Jordan's inequality and its applications to the improvement of the Yang Le inequality II},
 Appl. Math. Lett. {\bf 20}(5), 532-538 (2007)


\bibitem{18}
{\sc L. Zhu}, {\sc J.-K. Hua}: {\em  Sharpening the Becker-Stark inequalities},
 J. Inequal. Appl. 931275 (2010)

\bibitem{19}
{\sc L. Zhu}: {\em  Sharp Becker-Stark-type inequalities for Bessel functions},
 J. Inequal. Appl. 838740 (2010)

\bibitem{20}
{\sc L. Zhu}: {\em A refinement of the Becker-Stark inequalities},
Mat. Zametki {\bf 933}, 401-406 (2013)

\bibitem{21}
{\sc B. Male\v sevi\' c}, {\sc I. Jovovi\' c}, {\sc B. Banjac}:
{\em A proof of two conjectures of Chao-Ping Chen for inverse trigonometric functions},
J. Math. Inequal. {\bf 11} (1),  151-162 (2017)

\bibitem{22}
{\sc T. Lutovac}, {\sc B. Male\v sevi\' c}, {\sc M. Ra\v sajski}:
{\em A new method for proving some inequalities related to several special functions},
arXiv: 1802.02082 (2018)
\end{thebibliography}
\end{document}